\documentclass[10pt,a4paper]{article}

\usepackage{latexsym}
\usepackage{verbatim}
\usepackage{amsthm}
\usepackage{amsfonts}
\usepackage{amsmath}
\usepackage{amssymb}
\usepackage{amscd}

\numberwithin{equation}{section}

\newtheorem{prop}{Proposition}[section]
\newtheorem{theorem}[prop]{Theorem}

\newtheorem{corollary}[prop]{Corollary}
\newtheorem{remark}[prop]{Remark}

\def\begeq{\begin{equation}}
\def\endeq{\end{equation}}

\def\<{\langle}
\def\>{\rangle}
\def\({\left(}
\def\){\right)}

\def\Ric{{\rm Ric}}

\def\KC{{\rm KC}}

\def\p{\partial}

\begin{document}

\title{Ricci Lower Bound for K\"ahler-Ricci Flow}

\author{Zhou Zhang \\
Univeristy of Sydney} 
\date{}
\maketitle

\noindent {\bf Abstract:} in this note, we provide 
some general discussion on the Ricci lower bound 
along K\"ahler-Ricci flow with singularity over 
closed manifold. 

\section{Introduction and Set-up}

In this note, we study a topic regarding the behaviour 
of Ricci curvature for K\"ahler-Ricci flows over closed 
K\"ahler manifolds. This relates to the general brief 
(or conjecture) that Ricci curvature would have uniform 
lower bound towards (most) finite or infinite time 
singularities. A typical result is as follows, which is 
a more English version of Theorem \ref{th:finite}. 
\begin{theorem}
For any K\"ahler-Ricci flow with finite time singularity, 
if the global volume is not going to zero towards the 
time of singularity, then the Ricci curvature can not 
have a uniform lower bound. 
\end{theorem}

K\"ahler-Ricci flow is nothing but Ricci flow with 
initial metric being K\"ahler. For this K\"ahler 
condition of the initial metric, the closed smooth 
manifold, $X$, would be given a complex structure, 
which is fixed for the consideration. So we still 
call this complex manifold by $X$. The smooth flow 
metric would always be K\"ahler with respect to $X$, 
as first observed by R. Hamilton. We consider $dim_
\mathbb{C} X=n\geqslant 2$.  

Thus the standard form of K\"ahler-Ricci flow, a 
directly transformation from Ricci flow to a metric 
form flow, is as follows over $X\times [0, S)$ (for 
some $S\in (0, \infty]$)  
\begin{equation}
\label{eq:rfk}
\frac{\p \omega(s)}{\p s}=-2\Ric\(\omega(s)
\), ~~~~\omega(0)=\omega_0,
\end{equation}
where $\omega_0$ is the metric form for the initial 
K\"ahler metric. The key advantage in our study of 
this flow, comparing with many earlier works, is 
that we no longer force cohomology restriction to 
$[\omega_0]$ when setting up this flow (and the 
equivalent scalar metric potential flow appearing 
shortly). This allows more applications for this 
kind of geometric flow techniques and more importantly, 
makes it possible to analyze degenerate situation. 
This idea first appeared in \cite{tsu} and has been 
rigorized and generalized in the works of \cite{t-znote}, 
\cite{song-tian} and continuations.     

We often perform the following time-metric scaling 
for the above evolution equation, 
$$\omega(s)=e^t\widetilde\omega_t, ~~~~s=
\frac{e^t-1}{2},$$
and arrive at an equivalent version of K\"ahler-Ricci 
flow over $X\times [0, T)$ (for $T=\log(1+2S)\in (0, 
\infty]$), 
\begin{equation}
\label{eq:krf} 
\frac{\p \widetilde\omega_t}{\p t}=-\Ric
(\widetilde\omega_t)-\widetilde\omega_t, 
~~~~\widetilde\omega_0=\omega_0,
\end{equation}
which is somewhat more convenient as explained in the 
following.  

To begin with, let's point out that the time scaling 
makes sure that these two flows would both exist up 
to some finite times (with finite time singularity) 
or exist forever. In the finite time singularity case, 
the metric scaling is by uniformly control positive 
constants, and so one can say the equivalence is 
stronger. When they both exist forever, the metric 
scaling might have significant impact on the flow 
metric. For example, in the case of $c_1(X)=0$ as 
studied in \cite{cao}, $\omega(s)$ converges at time 
infinity to the unqiue Ricci-flat metric in the K\"ahler 
class $[\omega_0]$ with the lower $CY$ standing for 
the more popular name as Calabi-Yau metric, while 
$(X, \widetilde\omega_t)$ shrinks to a metric point 
at time infinity. Meanwhile, in other cases, $\omega
(s)$ has volume tending to infinty while $\widetilde
\omega_t$ has uniformly controlled volume. In principle, 
(\ref{eq:krf}) always has the cohomology information 
$[\widetilde\omega_t]$ under uniform control, as is 
clear in following discussion. 

{\it For the rest of this note, we focus on (\ref
{eq:krf}).}
 
Let's fix the convention that $[\Ric]=c_1(X)$, and 
then one can reduce (\ref{eq:krf}) to an ODE in the 
cohomology space $H^{1, 1}(X, \mathbb{R}):=H^2(X; 
\mathbb{R})\cap H^{1, 1}(X; \mathbb{C})$. It's easy 
to solve it and we end up with 
$$[\widetilde\omega_t]=-c_1(X)+e^{-t}\([\omega_0]+
c_1(X)\)$$
which stands for an interval in this vector space 
with two endpoints being $[\omega_0]$ for $t=0$ 
and $-c_1(X)$ formally for $t=\infty$. It's easy 
to see that for (\ref{eq:rfk}), $[\omega(s)]$ 
would evolve linearly in general, which might be a 
simpler function but the control is not as uniform.  

However, for either one, the optimal existence result 
in \cite{t-znote} tells that the flow metric exists 
as long as the class from the above consideration 
stays inside the open cone, consisting of all K\"ahler 
classes (called K\"ahler cone). In this note, we use 
$\KC(X)$ to denote the K\"ahler cone of $X$, and then 
we call its closure (in the finite dimensional vector 
space $H^{1, 1}(X, \mathbb{R})$), $\overline{\KC(X)}$, 
is sometime called the numerically effective cone, using 
a terminology borrowed from Agebraic Geometry. Now the 
optimal existence result simply says 
$$T=\sup\{t\,|\,-c_1(X)+e^{-t}\([\omega_0]+c_1(X)\)
\in \KC(X)\}$$
is the largest time to consider classic solution for 
(\ref{eq:krf}). This would be our definition for $T$ 
for the rest of this work, which takes value in $(0, 
\infty]$. Also from the result in \cite{t-znote}, 
we've already shown that singularity only happens when 
the class hits the boundary of this cone, at either 
finite or infinite time. In other words, $[\omega_T]
\in \overline{\KC(X)}\setminus \KC(X)$. Of course, 
now the interesting and in general hard question is 
to see how the non-K\"ahler feature of this class 
at the boundary of $\KC(X)$ and the behaviour of the 
K\"ahler-Ricci flow would interact with each other. 

The study of this topic, as well as for many other 
things regarding K\"ahler-Ricci flow, usually makes 
of the scalar version of the K\"ahler-Ricci flow. 
For (\ref{eq:krf}), we define the following background 
form,   
$$\omega_t=-\Ric(\omega_0)+e^{-t}\(\omega_0+\Ric
(\omega_0)\)$$
compatible with the notation $\omega_0$. The whole 
point is that $[\omega_t]=[\widetilde\omega_t]$, and 
so $\widetilde\omega_t=\omega_t+\sqrt{-1}\p\bar\p u$. 
It's not so hard to prove the following scalar 
evolution equation for metric potential $u$ over 
$X\times [0, T)$ is equivalent to (\ref{eq:krf}), 
\begin{equation}
\label{eq:skrf} 
\frac{\p u}{\p t}=\log\frac{\widetilde\omega^n_t}
{\omega^n_0}-u=\log\frac{(\omega_t+\sqrt{-1}\p\bar
\p u)^n}{\omega^n_0}-u, ~~~~u(\cdot, 0)=0.
\end{equation}
This evolution equation can be reformulated as 
\begin{equation}
\label{eq:cma} 
(\omega_t+\sqrt{-1}\p\bar\p u)^n=e^{\frac{\p u}
{\p t}+u}\omega^n_0,
\end{equation}
which is why sometimes we also call it the complex 
Monge-Amp\`ere equation type of K\"ahler-Ricci flow.

\vspace{0.1in}

Now let's also clarify that the uniform Ricci lower 
bound mentioned at the beginning. It means that there 
exists some constant $C$ such that
$$\Ric(\widetilde\omega_t)\geqslant -C\widetilde
\omega_t$$
uniformly for $t\in [0, T)$. In general, Ricci 
curvature being bounded from below provides some 
control on the metric and topology of the underlying 
manifold. So such a control is certainly favourable 
as an assumption and interesting as a result in the 
study of Ricci flow (see in \cite{bing} for examples 
of such assumptions).    

\vspace{0.1in}

\noindent{\bf Note:} in the following, $C$ always 
stands for a (positive) constant, which might be 
different at places. 

\section{Finite Time Singularity} 

In this section, we consider the case of $T<\infty$. 
Then clearly, $[\omega_T]\in\overline{\KC(X)}\setminus 
\KC(X)$. The following is the main result. 
\begin{theorem}
\label{th:finite} 
Consider (\ref{eq:krf}) with finite time singularity, 
i.e., $T<\infty$. If $[\omega_T]^n>0$, then the Ricci 
curvature can NOT have a uniform lower bound, i.e., 
there is NO constant $D>0$ such that $\Ric(\widetilde
\omega_t)\geqslant -D\widetilde\omega_t$ uniformly 
for $t\in [0, T)$. 
\end{theorem}
  
The proof is a combination of techniques from \cite
{r-blow-up} and \cite{weak-limit}. We begin with 
some general situation and finally specify to the 
case of the above theorem to prove it.  

The following is observed earlier as in Remark 2.3 of 
\cite{weak-limit}. It's clear that $[\omega_t]^n=
[\widetilde\omega_t]^n>0$ for $t\in [0, T)$, and we 
also have $[\widetilde\omega]^n=[\omega_t]^n\to[\omega
_T]^n$ as $t\to T$. So $[\omega_T]^n\geqslant 0$. In 
exactly the same manner, we see $[\omega_T]^{n-k}\cdot
[\omega_0]^k\geqslant 0$ for $k=1, \cdots, n-1$. Now 
rewrite $\omega_t$ as follows, 
$$\omega_t=\(\frac{1-e^{-t}}{1-e^{-T}}\)\omega_T+\(
\frac{e^{-t}-e^{-T}}{1-e^{-T}}\)\omega_0$$
and it is then obvious that for $t\in [0, T]$, 
$$[\omega_t]^n\thicksim (T-t)^K$$
where $K$ is defined as follows, 
$$n\geqslant K:=\min\{k\in\{0, 1, 2, \cdots, n\}\,|\, 
[\omega_T]^{n-k}\cdot[\omega_0]^k>0\},$$
which is well-defined since $[\omega_0]^n>0$. Here 
$A\thicksim B$ for $A$ and $B$ both non-negative 
means $\frac{1}{C}B\leqslant A\leqslant CB$ for 
some positive constant $C$.

\vspace{0.1in}

\noindent{\bf Note:} when $[\omega_T]^n>0$, $K=0$ 
and $\frac{1}{C}\leqslant [\omega_t]^n\leqslant C$ 
for constant $C>0$.

\vspace{0.1in}

As in Subsection 2.3 of \cite{weak-limit}, after 
assuming $\Ric(\widetilde\omega_t)\geqslant -C
\widetilde\omega_t$ for some constant $C>0$, 
plugging it into (\ref{eq:krf}) gives $\frac{\p 
\widetilde\omega_t}{\p t}\leqslant C\widetilde\omega$. 
Then noticing $T<\infty$, we arrive at $\widetilde
\omega_t\leqslant C\omega_0$ and so $\Ric(\widetilde
\omega_t)\geqslant -C\omega_0$.

The equivalent equations (\ref{eq:krf}) and (\ref
{eq:skrf}) give 
$$\Ric\(\widetilde\omega_t\)=-\frac{\p \widetilde\omega_
t}{\p t}-\widetilde\omega_t=\Ric(\omega_0)-\sqrt{-1}
\p\bar\p \(\frac{\p u}{\p t}+u\),$$
and so one arrives at
$$C\omega_0+\sqrt{-1}\p\bar\p\(-\frac{\p u}{\p t}-u\)
\geqslant 0.$$
Thus we can apply the classic result in \cite{tian-87} 
and have a constant $\alpha>0$ depending only on $(X, 
\omega_0)$ such that for $t\in [0, T)$, 
$$\int_X e^{\alpha\(\sup_X(-\frac{\p u}{\p t}-u)+
(\frac{\p u}{\p t}+u)\)}\omega^n_0\leqslant C.$$
Of course, we could make sure $\alpha\leqslant 1$. 
This gives
$$\inf_X\(\frac{\p u}{\p t}+u\)\geqslant\frac{1}
{\alpha}\log\(\frac{1}{C}\int_X e^{\alpha(\frac
{\p u}{\p t}+u)}\omega^n_0\).$$

As summarized in \cite{weak-limit}, we have $\frac
{\p u}{\p t}+u\leqslant C$, and so  
\begin{equation}
\begin{split}
\int_X e^{\alpha(\frac{\p u}{\p t}+u)}\omega^n_0
&= e^{\alpha C}\int_X e^{\alpha(\frac{\p u}
{\p t}+u-C)}\omega^n_0 \\
&\geqslant e^{\alpha C}\int_X e^{\frac{\p u}{\p t}+
u-C}\omega^n_0 \\
&\geqslant C\int_X e^{\frac{\p u}{\p t}+u}\omega^n_
0 \\
&= C[\widetilde\omega_t]^n=C[\omega_t]^n\geqslant 
C(T-t)^K \nonumber
\end{split}
\end{equation}
where $\alpha\leqslant 1$ is applied for the second 
step. So we conclude that for $t\in [0, T)$, 
$$\inf_X\(\frac{\p u}{\p t}+u\)\geqslant -C+\frac
{K}{\alpha}\log(T-t)$$
and so 
$$\frac{\p u}{\p t}+u\geqslant -C+\frac{K}{\alpha}
\log(T-t)$$
for $\alpha\in (0, 1]$ depending only on $(X, 
\omega_0)$. Directly applying Maximum Principle 
to (\ref{eq:skrf}), we have $u\leqslant C$ and so
\begin{equation} 
\label{ieq:volume-lower}
\frac{\p u}{\p t}\geqslant -C+\frac{K}{\alpha}
\log(T-t).
\end{equation}
So in a way slightly different from that in \cite
{weak-limit}, we have $u\geqslant -C$.  

The above estimate provides a pointwise lower 
bound of the volume form, $\widetilde\omega^
n_t=e^{\frac{\p u}{\p t}+u}\omega^n_0$. 
Combining with the metric upper bound, we arrive 
at the following proposition. 
\begin{prop}
\label{prop:metric-finite}
Consider (\ref{eq:krf}) with singularity at some 
finite time $T$. If $\Ric(\widetilde\omega_t)
\geqslant -D\widetilde\omega_t$ for some constant 
$D>0$ and $t\in [0, T)$, then 
$$(T-t)^\beta\omega_0\leqslant \widetilde\omega_t
\leqslant C\omega_0$$
for positive constants $\beta$ and $C$ depending 
on $X$, $\omega_0$, $T$ and $D$.
\end{prop}
  
Now we restrict to the case of Theorem 
{\ref{th:finite}}. In this case, $K=0$ 
and so the above lower bound for $\frac
{\p u}{\p t}$, (\ref{ieq:volume-lower}) 
is uniform. Hence, the metric control in 
Proposition {\ref{prop:metric-finite}} 
is also uniform. The argument in \cite
{r-blow-up} can be used to draw contradiction 
with the existence of finite time singularity 
at $T$. Theorem {\ref{th:finite}} is 
thus proven.  

\begin{remark}
This theorem indicates that for the problem in 
\cite{weak-limit} on general weak limit, when 
Ricci curvature has uniform lower bound, the 
discussion there is actually only for the global 
volume collapsed case. This again stresses the 
point that the discussion of collapsed case is 
the core for the topic of weak limit in general.   

It's worth pointing out that there are numerous 
examples satisfying the assumption of this theorem. 
For instance, we have the case discussed in \cite
{s-w}. In fact, such manifold of $X$ belongs to 
the class of so-called manifolds of general type, 
indicating ``majority". 

The problem on finite time singularity has been 
studied extensively since R. Hamilton's original 
work \cite{ham4}. The Ricci lower bound assumption 
and N. Sesum's result on the blow-up of Ricci 
curvature for finite time singularity of Ricci 
flows over closed manfolds in \cite{sesum} 
automatically gives the blow-up of the scalar 
curvature, which is conjectured in general and 
proven for K\"ahler case in \cite{r-blow-up}. 
Our theorem here actually shows that the 
situation has to be more complicated, at least 
in the global volume non-collapsed case.  
\end{remark}

\section{Infinite Time Singularity}

Now we consider the infinite time singularity case, 
i.e., $T=\infty$ and $[-\Ric(\omega_0)]=-c_1(X)\in
\overline{\KC(X)}\setminus \KC(X)$. When $X$ is 
projective, it is then a minimal manifold. Again 
we assume $\Ric(\widetilde\omega_t)\geqslant -D
\widetilde\omega_t$, which clearly gets weaker as 
$D$ gets larger. Our discussion below is separated 
into cases for increasing $D$ value, with the 
conclusion getting weaker. 

\begin{itemize}

\item $D<1$

In this case, (\ref{eq:krf}) gives $\frac{\p \widetilde
\omega_t}{\p t}\leqslant (D-1)\widetilde\omega_t$, and 
so $\widetilde\omega_t\leqslant e^{(D-1)t}\omega_0$.

So clearly $-c_1(X)=[-\Ric(\omega_0)]=0$, and the result 
in \cite{cao} can be scaled to provide a very satisfying 
description of the flow metric as follows. Using the 
notations in Section 1, $\omega(s)=e^t\widetilde\omega_t$ 
converges exponentially fast (for example, Section 9.3 in 
\cite{thesis}) to the Ricci-flat K\"ahler metric $\omega_
{CY}$, where this exponentially fast convergence is with 
respect to the parameter $s=\frac{e^t-1}{2}$. So as smooth 
forms, $\Ric(\widetilde\omega_t)$ is $e^{-s}$-small while 
$\widetilde\omega_t$ is $e^{-t}$-positive, and so the above 
Ricci lower bound is obviously true for large time. 

\item $D=1$

In this case, (\ref{eq:krf}) gives $\frac{\p \widetilde
\omega_t}{\p t}\leqslant 0$, and so $\widetilde\omega_t
\leqslant\omega_0$.

$\Ric(\widetilde\omega_t)+\widetilde\omega_t\geqslant 0$ 
also tells us that the corresponding cohomology class 
$$-c_1(X)+\(-c_1(X)+e^{-t}\bigl([\omega_0]+c_1(X)\bigr)
\)=e^{-t}([\omega_0]+c_1(X))\in\overline\KC(X),$$ 
and so $[\omega_0]+c_1(X)\in\overline\KC(X)$, providing 
a topological restriction.   

The above uniform metric upper bound allows most of 
the discussion in Section 2 to be carried through. 
Together with $\Ric(\widetilde\omega_t)\geqslant-
\widetilde\omega\geqslant -\omega_0$,  
$$\Ric\(\widetilde\omega_t\)=-\frac{\p \widetilde\omega_
t}{\p t}-\widetilde\omega_t=\Ric(\omega_0)-\sqrt{-1}
\p\bar\p \(\frac{\p u}{\p t}+u\),$$
will give us   
$$C\omega_0+\sqrt{-1}\p\bar\p\(-\frac{\p u}{\p t}-u\)
\geqslant 0.$$
Again apply the classic result in \cite{tian-87} to get 
constant $\alpha>0$ depending only on $(X, \omega_0)$ 
such that for $t\in [0, \infty)$, 
$$\int_X e^{\alpha\(\sup_X(-\frac{\p u}{\p t}-u)+
(\frac{\p u}{\p t}+u)\)}\omega^n_0\leqslant C.$$
Of course, we could make sure $\alpha\leqslant 1$. 
This gives
$$\inf_X\(\frac{\p u}{\p t}+u\)\geqslant\frac{1}
{\alpha}\log\(\frac{1}{C}\int_X e^{\alpha(\frac
{\p u}{\p t}+u)}\omega^n_0\).$$

As summarized in \cite{weak-limit}, we still have 
$\frac{\p u}{\p t}\leqslant C$ and $u\leqslant C$, 
and so in the same way as in Section 2, we arrive 
at  
$$\int_X e^{\alpha(\frac{\p u}{\p t}+u)}\omega^n_0
\geqslant C[\omega_t]^n.$$ 

Repeating the same discussion at the beginning of 
Section 2, we have $[-\Ric(\omega_0)]^{n-k}
\cdot [\omega_0]^k\geqslant 0$ for $k\in\{0, 1, 
\cdots, n\}$, where $-\Ric(\omega_0)$ can be viewed 
as $\omega_T$ for $T=\infty$. Furthermore, $[\omega_
t]^n\thicksim e^{-Kt}$ with   
$$n\geqslant K:=\min\{k\in\{0, 1, 2, \cdots, n\}\,|
\, [-\Ric(\omega_0)]^{n-k}\cdot[\omega_0]^k>0\},$$
which is well-defined since $[\omega_0]^n>0$. So we 
conclude that for $t\in [0, \infty)$, 
$$\inf_X\(\frac{\p u}{\p t}+u\)\geqslant -\frac{K}
{\alpha}t-C,$$ 
and so 
$$\frac{\p u}{\p t}+u\geqslant -\frac{K}{\alpha}t-
C$$ 
for $\alpha\in (0, 1]$ depending only on $(X, \omega_
0)$. This provides a pointwise lower bound of the 
volume form, $\widetilde\omega^n_t=e^{\frac{\p u}
{\p t}+u}\omega^n_0$. Combining with the metric 
upper bound, we arrive at the following proposition. 
\begin{prop}
\label{prop:metric-infinite}
Consider (\ref{eq:krf}) with the solution existing 
forever but having infinite time singularitiy. If 
$\Ric(\widetilde\omega_t)\geqslant -\widetilde\omega_
t$ for $t\in [0, \infty)$, then $[\omega_0]+c_1(X)
\in\overline\KC(X)$ and 
$$e^{-\beta t}\omega_0\leqslant \widetilde\omega_t
\leqslant \omega_0$$
for some positive constant $\beta$ depending on $(X, 
\omega_0)$.
\end{prop}
  
If we further assume $K=0$, i.e., $[-\Ric(\omega_0)
]^n>0$, then it's the global volume non-collapsed 
case and the metric bound from the above proposition 
is uniform. As in \cite{r-blow-up}, this implies 
$[-\Ric(\omega_0)]=-c_1(X)\in \KC(X)$, which 
contradicts the infinite time sinuglarity assumption 
which indicates $[-\Ric(\omega_0)]\in \overline{\KC
(X)}\setminus \KC(X)$. Let's summarize it in the 
following corollary, which is similar to Theorem 
{\ref{th:finite}} but not as neat. 
\begin{corollary}
\label{infinite}
Consider (\ref{eq:krf}) with the solution exists 
forever but having infinite time singularity. If 
$\Ric(\widetilde\omega_t)\geqslant -\widetilde
\omega_t$ for $t\in [0, \infty)$, then $c_1(X)^n
=0$.  
\end{corollary}

\item $D>1$

This is the general case. We can only have 
$\widetilde\omega_t\leqslant e^{Ct}\omega_0$ 
for some $C>0$. Then 
$$-Ce^{Ct}\omega_0\leqslant \Ric(\widetilde
\omega_t)=\Ric(\omega_0)-\sqrt{-1}\p\bar\p\(
\frac{\p u}{\p t}+u\),$$
and we could only have 
$$C\omega_0-\sqrt{-1}\p\bar\p \(e^{-Ct}\(
\frac{\p u}{\p t}+u\)\)\geqslant 0.$$
Applying the same discussion as before for 
$e^{-Ct}(\frac{\p u}{\p t}+u)$ only gives 
$$\frac{\p u}{\p t}+u\geqslant -Ce^{Ct}.$$
Hence, the metric bound corresponding to 
Proposition \ref{prop:metric-infinite} is
$$e^{-Ce^{Ct}}\omega_0\leqslant\widetilde
\omega_t\leqslant e^{Ct}\omega_0,$$
which is not enough to draw any decent 
conclusion. 

\begin{remark}
With this general lower bound of Ricci curvature 
for infinite time singularity case, when $X$ is 
a projective manifold of general type, i.e., 
$(-c_1(X))^n>0$, by the results in \cite{bound-r}, 
one has the Ricci curvature being bounded from 
both sides and $\frac{\p u}{\p t}+u\geqslant -C$. 
Thus the metric bound can be improved to $Ce^{-Ct}
\omega_0\leqslant\widetilde\omega_t\leqslant e^{Ct}
\omega_0$, not yet good enough. Notice that by 
Corollary \ref{infinite}, $D$ has to be strictly 
bigger than one to have a reasonable Ricci lower 
bound assumption for this case. 
\end{remark}

\end{itemize}

\vspace{0.2in}

\noindent {\it Email: zhangou@maths.usyd.edu.au}

\end{document}